\documentclass[12pt,a4paper]{article}
\usepackage[curve]{xypic}
 \usepackage[usenames,dvipsnames]{pstricks}
\usepackage{epsfig}
\usepackage[active]{srcltx}

\usepackage{amsthm}
\usepackage{amsfonts}
\usepackage{amscd}
\usepackage{amsmath}
\usepackage{amssymb}
\usepackage{color}
\usepackage{xcolor}
\usepackage{array}
\usepackage{epsfig}
\usepackage{graphics}
\usepackage{graphicx}
\usepackage{verbatim}
\usepackage{frcursive}
\usepackage[T1]{fontenc}
\usepackage{mathrsfs}

\usepackage{geometry}

\def\bfone{{\bf 1}}

\def\vecb0{\vec{\bf 0}}

\def\C{\mbox{$\mathbb C$}}

\def\nc{\newcommand}
\nc{\mystr}{{\rule[-1.5ex]{0ex}{4.5ex}{}}}
\def\mystrut{{\rule[-2ex]{0ex}{4.5ex}{}}}

\def\ld{ \left|\,\begin{matrix} }
 \def\rd{\end{matrix}\,\right|}

\def\rr{\rightarrow}

\def\dMT {d_{\rm MT}}
\def\dR {d_{\rm R}}
\def\SBV {{\rm SBV}}

\newcommand{\calS}{\mathcal S}

\def\rhosp{\rho_{sp}}
\def\rhoinf{\rho_{\infty}}

 \def\lv{ \left(\begin{matrix} }
 \def\rv{\end{matrix}\right)}

\def\ds{\displaystyle }

\newtheorem{newthm}{Theorem}
\newtheorem{theorem}{Theorem}[section]
\newtheorem{definition}{Definition}[section]
\newtheorem{lemma}[theorem]{Lemma}
\newtheorem{proposition}[theorem]{Proposition}
\newtheorem{corollary}[theorem]{Corollary}
\newtheorem{remarks}[theorem]{Remarks}
\newtheorem{example}[theorem]{Examples}

\newcommand{\REFEQN}[1] { \begin{equation}\label{#1} }
\newcommand{\ENDEQN}{\end{equation}}
\newcommand{\REFTHM}[1] { \begin{theorem}\label{#1} }
\newcommand{\ENDTHM}{\end{theorem}}
\newcommand{\REFDEF}[1] { \begin{definition}\label{#1} }
\newcommand{\ENDDEF}{\end{definition}}
\newcommand{\REFNTH}[1] { \begin{newthm}\label{#1} }
\newcommand{\ENDNTH}{\end{newthm}}
\newcommand{\REFPROP}[1]{\begin{proposition}\label{#1} }
\newcommand{\ENDPROP}{\end{proposition} }
\newcommand{\REFLEM}[1]{\begin{lemma}\label{#1} }
\newcommand{\ENDLEM}{\end{lemma} }
\newcommand{\REFCOR}[1]{\begin{corollary}\label{#1} }
\newcommand{\ENDCOR}{\end{corollary} }

\def\whA{\widehat{A}}
\def\whX{\widehat{X}}
\def\whL{\widehat{L}}
\def\whLk{\widehat{L_k}}

\def\whx{\widehat{x}} 
 \def\whf{\widehat{f}} \def\whI{\widehat{I}}
\def\whc{\widehat{c}} \def\whu{{\widehat{u}}} \def\whv{{\widehat{v}}}

\def\htop{h_{\rm top}}

\def\+{\text{\rm\tiny +}}   \def\-{\text{\small -}}

\DeclareMathAlphabet{\mathantt}{OT1}{antt}{m}{it}
\DeclareMathAlphabet{\mathpzc}{OT1}{pzc}{m}{it}

 \nc{\myhh}{\mbox{\cursive h}}
 \nc{\myh}{ \ell }
 \nc{\whh}{ \widehat{\ell} }
\nc{\myM}{\mbox{\cursive M}}
\nc{\normbv}[1]{\| #1 \|_{\rm BV}}

\def\SBV{ {\rm SBV}}

\def\myH{\myhh}
\nc{\Tr}{\mbox{\rm Tr}}
\nc{\tr}{\mbox{\rm tr}}
\nc{\mytrace}[1]{\Tr \left( \sstrut #1 \right)}

\nc{\mat}[2]  {\left(  \! \begin{array}{#1} #2 \end{array}\! \right)}

\DeclareMathOperator*{\var}{var}

\def\mystrut  {\rule[0ex]{0ex}{1.55ex}{}}
\nc{\myvar}[1]{\var_{\mystrut #1}}
\nc{\mychi}[1]{\chi_{\mystrut #1}}
\def\ds{\displaystyle}

\nc{\mydel}[1]{\Delta_{\mystrut #1}}
\nc{\mysig}[1]{\sigma_{\mystrut #1}}
\nc{\dual}[1]{\langle #1 \rangle}

\def\sigsp{\sigma_{\rm sp}}

\begin{document}
\title{The Milnor-Thurston
determinant and the
 Ruelle transfer operator.}
\author{  H.H. Rugh }
\maketitle

\begin{abstract} 
 The topological entropy $\htop$ of a continuous
 piecewise monotone interval map
measures the exponential growth in the number of
monotonicity intervals for iterates
of the map. 
 Milnor and Thurston showed that $\exp(-\htop)$
 is the smallest zero 
of an analytic function, now coined the Milnor-Thurston determinant,
 that keeps track of relative positions
of forward orbits of critical points.
On the other hand $\exp(\htop)$ equals the spectral radius of a 
Ruelle transfer operator $L$, associated with the map.
Iterates of $L$ keep track of inverse orbits of the map.
For no obvious reason, a Fredholm determinant for
the transfer operator has not only the same leading zero  as the
M-T determinant but all peripheral (those lying in the unit disk)
zeros are the same.

The purpose of this  note is to show that
on a suitable function space, the
dual of the Ruelle transfer operator has a 
regularized determinant,
identical to the Milnor-Thurston determinant,
hereby providing a natural explanation 
for the above puzzle.
\end{abstract}

\section{Introduction}
In the 70s, Milnor and Thurston  came up with an intriguing
way of computing the topological entropy $\htop$ 
for a continuous piecewise monotone interval map $f$.
They invented a "kneading-matrix", $M(t)$,
a finite  matrix-valued powerseries 
in an auxiliary variable $t$.
The matrix keeps track of 
relative positions of forward orbits of critical points relative
to the critical points themselves. The "Milnor-Thurston" determinant,
$\dMT(t)=\det M(t)$,
defines an analytic function
in the unit disk and they showed that if $\htop>0$ then
$t_*=\exp(-\htop) \in (0,1)$ is a zero of $\dMT(t)$.
The zero is `extremal' in the sense that it is the smallest in absolute
value.
Their result is computationally useful,
as the effort in a direct computation of the number of monotonicity
intervals of $f^n$ grows 
 like $e^{n \htop}$ whereas computing $M(t)$ to a similar accuracy
only grows linearly with $n$.

The above-mentioned extremal property
was shown indirectly by proving
that $\dMT(t)$ is identical to a certain  Lefschetz zeta-function,
generated from the periodic orbits for $f$.
The identity is shown in one (simple) case and
using continuous deformations of the map it is shown 
that the identity persists under such deformations.
This part is rather puzzling, as no arguments indicate
why one would expect such a relation. 

Also in the 70s, Ruelle 
used transfer operators to define 
(generalized) Fredholm determinants and
zeta-functions. In our context, 
an operator $L$ would act upon functions of bounded variation.
Iterates
of $L$ keep track
of inverse orbits of points and
its spectral radius equals precisely $1/t_*=\exp(\htop)$.
Baladi-Keller showed \cite{BK}, that
an associated determinant
$\dR(t) = \det(1-tL)$ is analytic in the unit disk and 
(at least for expanding maps) it is easy to see that
the zeros within the unit disk are the same as  for
the Lefschetz zeta function whence also for
$\dMT(t)$. This suggests  a deep relation
between the 
the M-T determinant and the
Ruelle-determinant of $L$, whence with the operator $L$ itself.

A step towards clarifying this relationship was suggested 
through a study of the dual operator
by Baladi and Ruelle \cite{BR}, but the functional
setup  makes explicit computions hard and requires operator-weights to be
globally continuous.
Gou\"ezel \cite{G},  elaborating
the functional setup further, 
 still requires weights to be continuous at periodic points.
Both of these results exclude partially the case studied by 
Milnor and Thurston.

The goal here
 is to show that the
determinant of Milnor-Thurston may be obtained by
restricting $L$ to a suitable small subspace $X$
of Bounded Variation functions, and calculate
an elementary regularized determinant of the dual operator $L'$.
We first exhibit a natural isometry between the dual space $X'$ and
a space $\whX$ of uniformly bounded functions of "point-germs".
This gives an explicit representation 
for the dual operator,
$\whL=S-PS$,
in which $S$ has spectral radius $1$ and
$P$ is of finite rank. When $|\lambda|>1$ (we call such a value
a 'peripheral' value) then $\lambda-S$ is
invertible through a von Neumann series. The elementary formula,
\[  \lambda- (S - PS) =
\left(1 + PS   \left(\lambda-S\right)^{-1} \right) 
   (\lambda-S) ,\]
shows that $\lambda$ is a spectral value iff
the operator 
$ (1 + PS(\lambda-S)^{-1} ) $ 
is non-invertible. 
As $P$ is of finite rank,  this
happens iff the following finite dimensional  determinant vanishes:
\begin{equation}
 \det \left( 1 + PS \left( \lambda - S\right)^{-1} \right) :=
 \det_{{\rm im}P} \left( P \sum_{k\geq 0}
        \lambda^{-k}   {S^k}  \right) .
  \end{equation}
Finally, setting
$\lambda=1/t$ this  turns out to be
the  Milnor-Thurston determinant.
An advantage of the present approach is that we may
exploit positivity of the Ruelle operator
(see Theorem \ref{thm positivity} below)
to bypass
the use of the Lifschetz zeta-function in the proof of the extremal property
of the zero.

%

\section{Step functions of  bounded variation.} 

Fix $-\infty\leq a<b\leq +\infty$.
The variation of a function $\phi : (a,b) \rr \C$
is:
\begin{equation}
   \var   \phi = \sup \left\{ \sum_{i} 
      \left| \phi(x_i)-\phi(x_{i+1}) \right|
       : a < x_1 < \ldots < x_N < b \right\} \ .
         \label{def var}
\end{equation}
We say that $\phi$ is of bounded variation (BV) iff $\var \phi < +\infty$.
For a BV-function right and left limits always exist and
following usual conventions we write
$\phi(x^+)=\lim_{t\rr 0^+} \phi(x+t)$ and
$\phi(x^-)=\lim_{t\rr 0^+} \phi(x-t)$ defined for $x\in[a,b)$ and
$x\in (a,b]$, respectively.
We assume no {a priori}
relation between the values of $\phi$ at $x^-$,$x$ and $x^+$
  (see also Remark \ref{remark choice} below).
We define the 'boundary' value of $\phi$ to be
\begin{equation}
   \partial \phi = \phi(a^+)+\phi(b^-)  \ \ .
 \end{equation}

In the following we write
$\whx=x^\epsilon$ for either  $x^+$ or $x^-$ and
call it a point germ with base point
$x$ and direction
 $\epsilon=+1$ or $\epsilon=-1$.
We order
point germs intertwined with 
base points by declaring that when $x<y$ then
      $x < x^+ < y^- < y$.
Any real segment may then be specified through
two point-germs by setting:
    $\dual{ \whu, \whv } = \left\{  x\in(a,b) : \whu < x < \whv \right\}$.
So 
e.g.\ $\dual{ u^-,v^-} = [u,v)$ (for $u<v$),
$\dual{ u^-, u^+} = \{ u \}$ and $\dual{\whu,\whu}=\emptyset$.



  \begin{definition}
  A map $\phi:(a,b) \rr \C$ is said to be a simple
  step function if there is a finite set $C=C_\phi\subset (a,b)$
  so that $\phi$ is constant on each connented component
  of $(a,b)\setminus C$.
   We write $\calS=\calS(a,b)$ for the space
  of simple step functions on $(a,b)$.
  \end{definition}

When $a<x_0<x<x_N<b$ then 
$2|\phi(x)|\leq 
\var \phi + |\phi(x_0)+\phi(x_N)|$. Letting
  $x_0\rr a^+$ and $x_N\rr b^-$
yields 
$\sup |\phi| \leq  \frac{1}{2} (\var{\phi} + |\partial \phi|)$.
One easily obtains~:

\begin{proposition}
   The quantity $\normbv{\phi} = \var{\phi} + |\partial \phi|$ defines
   a norm on $S$. We let $X$ denote the completion of $S$ under this
   norm and call $X=\SBV(a,b)$ 
   the space of Step functions of Bounded Variation on $(a,b)$.
   One has 
   \ $\sup|\phi| \leq \frac12 \normbv{\phi}$
   \ for $\phi\in X$.
\end{proposition}



We define for every $a^+\leq \whu \leq b^-$
the following "base" step-function (which is everywhere non-zero):
  \begin{equation}
    \sigma_{\whu} (x) = \frac12 \times \left\{
	    \ds \begin{array} {ll}
	       + 1 \ \ & , \ \whu < x < b  \\
	       - 1 \ \  &, \ a < x < \whu 
	    \end{array} \right. 
  \end{equation}
In particular,
$\;\sigma_{a^+}(x) = - \; \sigma_{b^-}(x) = \frac12\;$
for all $a<x<b$. In terms of characteristic functions one has:
$\sigma_\whu (x) = \frac12 \left( 
- \mychi{\dual{a^+,\whu}} (x)
+ \mychi{\dual{\whu,b^-}} (x) \right)$.
	 

Let $X'$ be the Banach dual of $X$, i.e.\ the space of
linear maps  $\myh:X \rr \C$ for which 
    $\dual{\myh, \phi} \leq C \normbv{\phi}, \ \ \forall \phi \in X$
    and some constant $C<+\infty$.
The norm $\|\myh\|_{X'}$
of $\myh$ is defined as
the smallest such constant.
Acting with $\myh$ upon a base
step function $\mysig{\whu}$ 
we obtain a representation of $\myh$ as a
function on the set of point germs $\whI=[a^+,b^-]$. 
We call this the $\sigma$-transform 
of $h'\in X'$ and write
\begin{equation}
   \whh (\whu) := \dual{\myh, \mysig{\whu}} \  , \ \ \whu\in \whI.
\end{equation}

As is readily seen, $\mysig{\whu}$ 
has norm one, so 
$\| \myh \|_\infty \leq \|\myh \|_{X'}$.
Thus, $\myh$ is an element of the space
 $\whX=B([a^+,b^-])$,
 the bounded functions on $\whI$ equipped with the
uniform norm $\|\cdot\|_\infty$.
Since $\dual{\myh,\mysig{b^-}} =
-\dual{\myh,\mysig{a^+}}$, we have $\myh(a^+)+\myh(b^-)=0$
so $\myh$
in fact belongs to the closed subspace
$\whX_0 = \{ H \in \whX : H(a^+)+H(b^-)=0 \}$.
But more is true:

\begin{proposition}
\label{prop duality}
The $\sigma$-transform
$\ell\in X' \mapsto  \whh \in \whX_0$ 
is a Banach space isometry.
\end{proposition}
Proof:
	The collection $\{\mysig{\whu}: a^+\leq \whu< b^-\}$
forms a base for $S$, the space of simple step functions.
In fact, given any simple step function $\phi$
there is a finite number of point germs
$a^+ =\whu_0 < \cdots < \whu_N <b^-$ so that:
\begin{equation}
\phi = 
 w_0 \frac12 \bfone +
\sum_{i=1}^N w_i \sigma_{\whu_i} 
 = \sum_{i=0}^N w_i \sigma_{\whu_i} .
   \label{decomp}
   \end{equation}

We omit a term involving $b^-$ since 
$-\mysig{b^-}=\mysig{a^+}$ and the latter is already included in the sum.
The decomposition is unique if every $w_i\neq 0$, $i > 0$.
In terms of this decomposition:
\begin{equation}
  \var  \phi = \sum_{i=1}^N |w_i| \ \ \mbox{and} \ \ \partial \phi = w_0 \ .
\end{equation}
To see this, note that the variation is bounded by
the right hand sum and as is readily checked,
the supremum
 in the definition (\ref{def var}) may be realized
by choosing a mesh that verifies:
$a<x_1<\whu_1<x_2<\whu_2< \ldots < \whu_N < x_{N+1} < b$.
The norm of $\phi$  is 
 given by:
\begin{equation}
   \normbv{\phi} = \sum_{i=0}^N |w_i|.
   \label{phi norm}
   \end{equation}
Given a simple step function $\phi\in \calS$ 
it may be decomposed as in 
(\ref{decomp}) 
and using the expression (\ref{phi norm}) for the norm we see that:
\begin{equation}
   \dual{\myh,\phi} = \sum_{i=0}^N  w_i  \myh (\whu_i) 
	  \end{equation}
\begin{equation}
  \left| \dual{\myh,\phi} \right| = 
    \sum_{i=0}^N | w_i|   \ |\myh (\whu_i) |
          \leq \sum_{i=0}^N |w_i|  \ \| \whh \|_\infty = 
	  \normbv{\phi}   \| \whh \|_\infty .
	  \end{equation}
Now, this extends by continuity
to every $\phi$ in $X$, the completion of $\calS$.
Conversely, 
any function $H\in \whX_0$ gives rise to
a linear functional on $\calS$
of the same norm,
by setting 
$ \dual{\myh,\phi} := \sum_{i=0}^N  w_i  H (\whu_i)$
and extending by continuity $\ell$ becomes a linear functional
$\myh\in X'$ of the same norm.
\qed

\begin{remarks}
 The space $\SBV(a,b)$ is isomorphic to
   $\C \oplus \ell^1((a,b))$
  (summable functions on the uncountable set of point germs + the boundary value).
 Not surprisingly, $X'$ is isomorphic to 
   $\C \oplus \ell^\infty((a,b))$.
\end{remarks}

\begin{proposition}
\label{remark spectrum}
Given a bounded linear operator
$A\in L(X)$ we obtain
a representation $\whA \in L(\whX_0)$
 through:
 \begin{equation}
   \left( \whA\, \whh \right)  (\whu) := 
\dual{A' \myh, \mysig{\whu}} = \dual{\myh, A \mysig{\whu}},
 \ \ \ \whu \in [a^+,b^-].
  \label{eq repr}
  \end{equation}
  
The operators $A$ and $A'$ (whence $\whA$)
have the same spectrum
\cite[III, Thm 6.22]{Kato}.

$\lambda$ is an isolated eigenvalue of finite
albraic multiplicity for $A$ iff it is for $A'$ (whence for $\whA$).
In this case, algebraic (and geometric) multiplicites are the same
for $A$, $A'$ and $\whA$.
\cite[III, Remark 6.23]{Kato}.
\end{proposition}

In the following we shall study a family of operators
where the representation of the dual operator takes a
particular simple form.


%

\section{A Ruelle transfer operator}

Let $\emptyset \neq I_k=(a_k,b_k) \subset (a,b)$
(below, $k$ will serve as an index)
and let $f_k :I_k \rr (a,b)$ be a continuous,
strictly monotone map. We write
 $s_k=s(f_k)\in \{\pm 1\}$ for the sign of monotonocity.
Also let $g_k\in \C$
(the topological entropy corresponds to
the case $g_k\equiv 1$).
Given a triple $(I_k,f_k,g_k)$ we associate
a bounded linear operator,  $L_k$ acting upon
 $\phi \in X=\SBV(a,b)$ and
defined as follows:
\begin{equation}
   L_k \phi \left(y\right)
    =  \ g_k  \  \;  \phi\! \circ\! f_{k}^{-1} (y) \  \; \mychi{f_kI_k}(y) \  ,
       \ \ y\in (a,b).
\end{equation}
Acting with $L_k$
e.g.\ upon a characteristic function $\mychi{J}$
of $J\subset (a,b)$ simplifies to:
\begin{equation}
   L_k \mychi{J} \left(y\right) 
       =  \ g_k \  \mychi{f_k(I_k\cap J)}(y)  .
\end{equation}

\begin{figure}
\hspace{0.025\paperwidth}
\includegraphics[width=0.3\paperwidth]{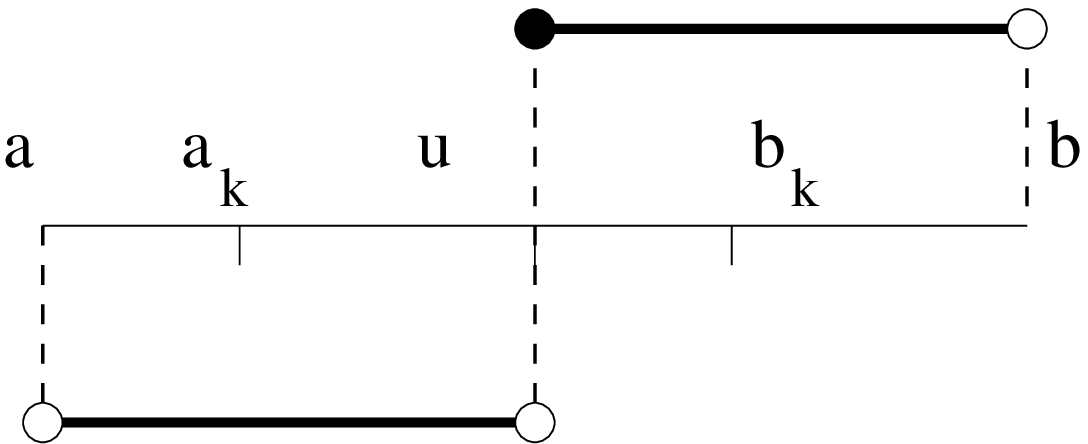}
\hspace{0.025\paperwidth}
\includegraphics[width=0.3\paperwidth]{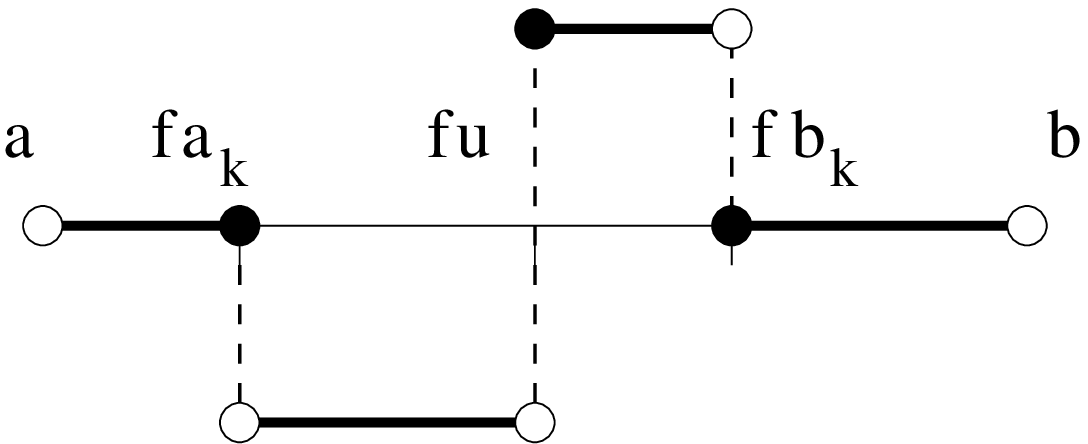}
\hspace{0.025\paperwidth}
\caption{Graph of $\sigma_{u^-}$  and
    $L_k \sigma_{u^-}$
for $a_k<u<b_k$ and 
    with $s_k=+1$.}
    \label{fig sigma}
\end{figure}

We extend  $f_k$ 
to a map of point germs
$\whu=u^\epsilon \in [a_k^+,b_k^-]$ by setting
$\whf_k(u^\epsilon) = (u')^{\epsilon'}$,
where 
$u'=  f_k(u^\epsilon)$
and $\epsilon'=s_k \epsilon$ (the direction of the
point germ is reversed precisely when $f_k$ reverses orientation).
By the results of the previous section,
the action of $L_k$ on an SBV-function is 
uniquely determined by its action upon 
a base  step function.
When $a_k <\whu < b_k$ the action 
of $L_k$ upon $\mysig{\whu}$ yields
a simple step function with
at most 3 jumps (see Figure \ref{fig sigma}) which may
therefore be written
as the sum of (at most) 3 base step functions. When $\whu$ is
outside of the interval $[a_k^+,b_k^-]$,
(at most) 2 such base step functions suffice.
Taking into account the reversal of orientation when $s_k=-1$ we
may summarize as follows~:

 \begin{equation}
 L_k  \mysig{\whu} = \frac12 \;s_k g_k \times
  \left\{ \begin{array}{lll}
        + \mysig{\whf_k a_k^+} -  \mysig{\whf_kb_k^-},
           & a \; < \whu < a_k\\[2mm]
      - \mysig{\whf_k a_k^+} -  \mysig{\whf_kb_k^-}
       +  2 \mysig{\whf_k \whu}
       \ \ \ 
      & a_k < \whu <b_{k}\\[2mm]
       -  \mysig{\whf_k a_k^+} +  \mysig{\whf_kb_k^-},
      & b_{k} < \whu < b
     \end{array} \right.  \ \ .
\end{equation}

The term  $\mysig{\whf_k \whu}$ only occurs
when $a_k<\whu<b_k$ while both
$\mysig{\whf_k a_k^+}$ and $\mysig{\whf_k  b_k^-}$
always appear, but with a sign depending upon the position of 
$\whu$ relative to $a_k$ and $b_k$, respectively.
We may collect terms and write:
\begin{equation}
  L_k \mysig{\whu} =
      s_k g_k \left[
      \mychi{ (a_k,b_k) }(\whu) \; \mysig{\whf_k\whu}
      -\mysig{a_k}(\whu)\;   \mysig{\whf_ka_k^+}
      +\mysig{b_k}(\whu)  \; \mysig{\whf_kb_k^-}
    \right] \ .
\end{equation}
Here, $\mysig{a_k}(\whu)=-\mysig{\whu}(a_k)$ takes the
value $+1/2$ when $\whu>a_k$ and $-1/2$ when $\whu<a_k$
(and similarly for $\mysig{b_k}(\whu)$).
The essential point is that the RHS is again a sum of elementary
step functions. Thus, acting with a linear functional on both sides 
we obtain a representation
for the dual of the  Ruelle operator:
\begin{eqnarray}
   \whLk \whh (\whu) &:=& \dual {L_k'\myh,\mysig{\whu}} 
      = \dual {\myh ,L_k\mysig{\whu}} \\
      & =& s_k g_k \left[
      \mychi{ (a_k,b_k) }(\whu) \whh(\whf_k\whu)
      -\mysig{a_k}(\whu)  \whh(\whf_ka_k^+)
      +\mysig{b_k}(\whu)  \whh(\whf_kb_k^-)
    \right] \ .
\end{eqnarray}
We extend this to any $H\in \whX$ by setting
   \begin{equation}
   \whLk H (\whu) = 
      s_k g_k \left[
      \mychi{ (a_k,b_k) }(\whu) H(\whf_k\whu)
      -\mysig{a_k}(\whu)  H(\whf_ka_k^+)
      +\mysig{b_k}(\whu)  H(\whf_kb_k^-)
    \right] \ .
\end{equation}
The constant function $\bfone$ is in the kernel of $\whLk$, since
for every $\whu\in[a^+,b^-]$:
\begin{equation}
 \whLk \bfone(\whu) = s_kg_k \left[ \mychi{ (a_k,b_k) }(\whu)
    - \mysig{a_k}(\whu)  +\mysig{b_k}(\whu) \right] \equiv  0.
       \label{kernel L}
\end{equation}

One also has
$\; \whLk H (a^+) = \frac12 s_k g_k
    \left( H(\whf_ka_k^+) - H(\whf_k b_k^-) \right) = - \whLk H( b^-)$.
Consequently,
$\;\whLk H (a^+)+\whLk H( b^-)=0$,
 so $\whLk$ not only preserves $\whX_0$ but maps $\whX$ into $\whX_0$.
In particular, any (non-zero) spectral property 
does not change by considering the extended operator.
In order to simplify expressions,
we introduce the following two operators acting upon $H\in \whX$:
\begin{eqnarray}
   S_k H(\whu) & =&  s_k g_k \mychi{ (a_k,b_k) }(\whu) H(\whf_k\whu), \\
   P_k H(\whu) & = &
         \mysig{a_k}(\whu) H(a_k^+) -\mysig{b_k}(\whu) H(b_k^-) .
\end{eqnarray}
In terms of these, our dual representation may be written as
\begin{equation}
    \whLk = S_k - P_k S_k .
    \end{equation}
The weighted
composition operator $S_k$ is (just) a bounded linear
weighted composition operator,
but the image of each $P_k$ is  spanned by the elements
$\mysig{a_k}$ and $\mysig{b_k}$
so $P_k$ is of rank two.
Our representation 
of the dual Ruelle operator is
 therefore a finite rank perturbation
of $S_k$.

As the reader may verify,
neither $S_k$ nor $P_k S_k$ need preserve $\whX_0$
(although their difference does).
It is, however, somewhat easier to work with
$S_k$ and $P_k$ acting upon $\whX$, rather than making 
identifications to restrict their action to $\whX_0$.

\subsection{A dynamical system}
We consider here the situation when $I_k=(a_k,b_k)=(c_k,c_{k+1})$ with
$a=c_0 < c_1 < \ldots < c_{d+1}=b$.
The intervals are disjoint and their union 
fill out $(a,b)$ apart from the cutting points $c_1,\ldots, c_d$.

We  associate as in the previous section 
a map, a weight and an operator to each interval.
We write $(I_k,f_k,g_k)$, $k=1,\ldots,d\;$ for the
collection of triples and denote by  $L=\sum_k L_k$ the
Ruelle transfer operator associated to this collection.

By our choice of cutting points,
we may extend the collection of $\whf_k$'s to a unique map $\whf$ 
defined upon the entire set of point germs $[a^+,b^-]$
(even if $f$ need not extend to a continuous map on $I$ !).
Signs and weights also extend to functions of $[a^+,b^-]$.
When $\whu\in \whI_k=[c_k^+,c_{k-1}^-]$.
	we write $s(\whu) = s_k$, $g(\whu)=g_k$  and
we have the expression for the
corresponding sum of composition operators:
\[ SH (\whu) = \sum_k S_k H (\whu)=
   \sum_k  s_k g_k \mychi{ (a_k,b_k) }(\whu)
      H \circ \whf_k  (\whu) = (s g)(\whu) H (\whf\; \whu).\]
Let us also write
$g^{(n)}(\whu) = g\circ \whf^{n-1}(\whu) \cdots g(\whu)$
(and similarly for $s^{(n)}$)
for the
product of weights (or signs)  along the orbit of $\whu$.
Then $ S^nH (\whu) = 
       (s g)^{(n)}(\whu) H (\whf^n \whu)$
and the spectral radius of $S$ is  given by
\[ \rho_\infty 
          = \lim_n \|S^n\|_\infty^{1/n}
          = \lim_n 
       \|g^{(n)} \|_\infty^{1/n}
       .\]
From the previous section we see that
the dual of the  Ruelle transfer operator, $L'=\sum_k L'_k$,
is a finite rank perturbation
of $S$ so the essential spectral radius of $L'$ 
(whence of $L$)
is not greater than $\rho_\infty$. 
If $\rho(L)>\rho_\infty$ then $L$ is quasi-compact
(\cite[VIII.8]{DS}) and any spectral value $\lambda$ with 
$|\lambda|>\rhoinf$ (we call such a value a peripheral spectral
value) must be an isolated eigenvalue of $L$ of finite 
algebraic multiplicity.


A function of the form $S_k H$  has support in $I_k$ so we have
$P_m S_k=0$ whenever $k\neq m$.
Whence, $\sum_k P_k S_k = P\; S$, where $P=\sum P_k$.
Furthermore,  we may recollect terms to get:
  \begin{align*}
  P H(\whu) 
    &= \sum_{k=0}^d \left[ \mysig{c_k}(\whu) H(c_k^+) -
        \mysig{c_{k+1}}(\whu) H(c_{k+1}^-) \right] 
	  =
      \mysig{a}(\whu) \Delta_{a}  H
      + 
      \sum_{j=1}^d  \mysig{c_j}(\whu) \Delta_{c_j} H \ ,
\end{align*} 
in which we have used $\mysig{b}(\whu)=-\mysig{a}(\whu)=1$ and defined
\begin{eqnarray}
     \Delta_{c_j} H &=& H(c_j^+) - H(c_j^-), \ \ j \geq 1 \  , 
        \label{Delta one} 
          \\
     \Delta_{a} H &=& H(a^+) + H(b^-) \ .
        \label{Delta two} 
     \end{eqnarray}
The image of $P$ is 
$d+1$ dimensional and as $\Delta_{c} \mysig{c'}=\delta_{c,c'}$,
 $P$  is also a  projection.
If $|\lambda|> \rhoinf$ 
then $\lambda-S$ is invertible through a von Neumann series so
$ \lambda-\whL = \left( 1 + PS (\lambda-S)^{-1} \right) (\lambda -S)$.
   Equivalently, writing $t=1/\lambda$:
 \begin{equation}
  1 -t \whL =  G(t)  (1  - t S),
   \label{WA formula}
   \end{equation}
where
\begin{equation}
G(t)
 = 1 + t \; PS (1-tS)^{-1}
 = (1-P) +  P (1-tS)^{-1}.
\end{equation}

We conclude that  $\lambda - \whL$ is non-invertible
(i.e.\ $\lambda=1/t$ is a spectral value)
iff the operator $G(t)$
is non-invertible. 

Since $(1-P)G(t)=(1-P)$ it suffices to look at the restriction of $G(t)$
to the image of $P$ which 
is spanned by
$\mysig{c_0}=\mysig{a}$ and
$\mysig{c_1},\ldots,\mysig{c_d}$.
Acting upon these $d+1$ functions
we get the following  $(d+1)\times(d+1)$ matrix representation
\[ G(t) \; \mysig{c_k} =  P (1-tS)^{-1} \mysig{c_k} =
     \sum_{j=0}^d \mysig{c_j} M_{jk}(t) \]
     where
 \begin{equation}
 M_{jk}(t)
   = 
      \Delta_{c_j} \left( \left( 1 - t  S\right)^{-1} 
         \mysig{c_k} \right)  = \Delta_{c_j} \Theta_{c_k}(\cdot,t)
	 \label{M matrix}
\end{equation}
and
\begin{equation}
   \Theta_{c_k}(\whx,t) = 
   \left( \left( 1-tS \right)^{-1} \mysig{c_k} \right) \left(\whx\right)
   =\sum_{t\geq 0} t^n (sg)^n(\whx) 
        \mysig{c_k}(\whf^n \whx),
    \label{kneading coor}
	\end{equation}
is equivalent to the so-called kneading coordinate of Milnor-Thurston.
Our matrix $M(t)$ is of 
size $(\ell+1_)\times(\ell+1)$ and is one version of the Milnor-Thurston 
kneading matrix. The original \cite{MT}
is, $\ell \times(\ell+1)$ dimensional) but
the associated M-T determinant
is identical to  the determinant of $M(t)$
as shown in  \cite[Appendix B]{RT}.

The value $\lambda=1/t$ is an eigenvalue 
iff $\dMT(t) := \det M(t) =0$.
Moreover, the order of the zero of the determinant equals
the algebraic multiplicity of the
 eigenvalue.
This last property is known as the Weinstein Aronszajn formula
(while not difficult to show it is somewhat lengthier,
so we refer to e.g.\ \cite[Ch. IV.6]{Kato} or 
\cite[Lemma 2.8]{R99} for a proof).
In summary we have shown:

\begin{theorem}
  Let $\sigsp(L)$ be the spectrum of the Ruelle transfer operator
  acting upon $X=\SBV(a,b)$. The peripheral spectrum,
  $\{z\in \sigsp(L): |z|>\rhoinf\}$ 
  consists of isolated eigenvalues of finite multiplicity only.
  Furthermore, $\lambda$ is a peripheral eigenvalue
  iff $t=1/\lambda$ is a zero of the Milnor Thurston
  determinant $\dMT(t)$ and the algebraic multiplicity of 
  $\lambda$ equals the order of the zero of $t$.
\end{theorem}

\begin{remarks}
  \label{remark choice}
  In the above we have used a representation of a simple step
  function $\phi$
in which the value at each point $x\in (a,b)$ is specified independently
of the limit-values at $x^+$ and $x^-$. Other possibilities are
to assign the value $\frac12\left( \phi(x^-)+\phi(x^-) \right)$
at $x$ or to declare two functions equivalent if they agree except
at a countable set.
The dual space of the completion in both of the latter cases
may be identified with  $\{H\in \ell^\infty([a,b]) : H(a)+H(b)=0\}$.
This is, however, 
inconvenient when the map is discontinuous and
the forward orbits of $c_j^+$ and
$c_j^-$ are distinct. With our approach the distinction is automatic
as it is build into the space of point-germs.
\end{remarks}

\subsection{Positivity and leading zero}

The most interesting situation is
when all weights are non-negative
and the Ruelle transfer operator  is
positive, i.e.\ preserves the cone of non-negative functions.
We may use this to give a direct simple proof of the
extremal properties
of the zero as alluded to in the introduction.
Let $Z_n$ denote the collection of monotonicity intervals
(also known as cylinder sets) of $f^n$.
We define the partition function  by
\begin{equation}
  \Omega_n= \sum_{\alpha \in Z_n} g^{(n)}_{|\alpha}
  \end{equation}
The particular case $g_k\equiv 1$ for every $k$, yields
the topological entropy through $\htop=
    \log \lim_n \Omega_n^{1/n}$.

\begin{theorem}
\label{thm positivity}
 Suppose that all $g_k\geq 0$. Then we have
 \begin{equation}
    \rho_1:= 
    \lim_n \Omega_n^{1/n} 
       = \lim_n \| L^n \bfone\|_\infty^{1/n}
       = \rho_{sp}(L). 
 \end{equation}
Moreover, if $\rho_{sp}(L)>\rho_\infty$ then 
  $\rho_{sp}(L)$ is an isolated eigenvalue of $L$ of finite multiplicity
  and $t^*=1/\rhosp(L)>0$ is the smallest zero (in absolute value) of
  the Milnor-Thurston determinant $\dMT(t)$.
\end{theorem}

Proof:
As $L^n \bfone(y) \leq \Omega_n$ for all $n$ and $y$ we have
 $\rho_\infty \leq \lim_n \|L^n \bfone\|^{1/n} \leq \rho_1$.

Consider the set $Z_n$ of (non-empty) monotonicity intervals for $f^n$.
Any open interval $\alpha\in Z_n$
may be written in terms of point germs as $\alpha=\dual{\whu,\whv}$ 
where each directed end point is either
a cutting point itself or it maps into 
a cutting point within at most $n-1$ iterations. 
Let $0\leq k<n$ be
the smallest such number for which $\whf^k (\whx) = \whc$ where
$\whc$ is either $c_j^-$ or $c_j^+$ for some $1\leq j\leq d$.
Now, the partition sum is the sum of $g^{(n)}_{|\alpha}$
over all $\alpha\in Z_n$.
We get twice the value if we instead sum
over all interval endpoints.
Using that $g^{(n)}(\whx) = g^{(k)}(\whx)  g^{(n-k)}(\whc_j)$ 
and rewriting the sum in view of the previous
remarks we get the identity (the boundary points $a^+$ and $b^-$
need a special treatment):
\begin{equation}
   2 \Omega_n = 
    \sum_{\dual{\whu,\whv} \in Z_n} \sum_{\whx\in \{\whu,\whv\}}
       g^{(n)}(\whx) = 
   g^{(n)}(a^+) + g^{(n)}(b^-) +
    \sum_{j=1}^d \sum_{k=0}^{n-1} 
        L^k \bfone  (c_j) 
       \sum_{\whc=c_j^{\pm}} g^{(n-k)} (\whc).
\end{equation}
We deduce the inequality
    $2\Omega_n \leq 2 \|g^{(n)}\|_\infty + 
        d \; \sum_{k=0}^{n-1} \|L^k \bfone\|_\infty \|g^{(n-k)}\|_\infty$
and then:
  \[ 
  \rho_1=\lim_{n\rr \infty} \Omega_n^{1/n} \leq
  \lim_{n\rr \infty} \|L^n \bfone\|_\infty^{1/n}
   \leq \|L^n \bfone\|_X^{1/n}\leq \rhosp(L). \]

If $\rhosp(L)=\rho_\infty$ we are through. So
suppose that $\rho_{sp}(L)>\rho_\infty$. As any peripheral spectral
value is en eigenvalue of finite multiplicity, there must be
$\lambda\in \C$ and a non-trivial $\phi\in X$ with $\lambda\phi=L\phi$
and $|\lambda|=\rhosp(L)$.
Taking absolute values we get by positivity of the operator
  \begin{equation}
|\lambda|\;  |\phi| =
     \left| \sum_k g_k \  \phi\circ f_k^{-1}
      \ \mychi{f_kI_k} \right| 
  \leq 
     \sum_k g_k \  |\phi|\circ f_k^{-1}
      \ \mychi{f_kI_k}  
   =    L |\phi| \ .
\end{equation}
Similarly $|\lambda|^n |\phi| \leq L^n |\phi| \leq L^n \bfone \|\phi\|_\infty$,
from which we deduce that 
\[
\rhosp(L)=|\lambda| \leq \lim_n \|L^n\bfone\|^{1/n}\leq \rho_1 \leq \rhosp(L).
\]

For the last assertion, consider the above $\phi$ and $\lambda$.
When $0<t<1/\rhosp(L)$ we get
$(\bfone -tL)^{-1} |\phi| =
\sum t^k L^k |\phi| 
  \geq \sum t^k |\lambda|^k |\phi|  = (1-t \; \rhosp(L))^{-1} |\phi|
  $ which diverges as $t\rr (1/\rhosp(L))^-$. This implies that
  $(\rhosp(L) \bfone - L)$ must be non-invertible and whence that
  $\rhosp(L)$ is an eigenvalue of $L$.
Then $t_*=1/\rhosp(L)$ is a zero of $\dMT(t)$,
and the smallest such since there can
be no spectral value larger than $\rhosp(L)$.\qed
\\

\begin{example}
 Consider the two triples of maps (with unit weights):

 $I_1=(0,\frac12), f_1(x)=2x, g_1=1$ and
 $I_2=(\frac12,1), f_2(x)=x-\frac12, g_2=1$.\\
 $\whf$ maps the four directed boundary points
 $0^+, \frac12^-, \frac12^+$ and $1^-$ to
 $0^+,1^-,0^+$ and $\frac12^-$, respectively.
      $\Theta_\frac12( \frac12^+,t)$ is calculated by
      looking at the
      forward orbit of $\frac12^+$ relative to $\frac12$.
   For example,
   as $\frac12^+>\frac12$ (which yields a plus sign)
   and $\whf(\frac12^+)=0^+=\whf(0^+) < \frac12$ (so a minus sign)
   we get 
   \begin{equation}
      \Theta_\frac12( \frac12^+,t)=
       \frac12 \left(1 - t - t^2 - t^3 -\cdots \right) =
       \frac12 \frac{1-2t}{1-t}. 
       \end{equation}
 Our $2\times 2$ kneading matrix, cf.\ 
  (\ref{Delta one}), (\ref{Delta two}) and
  (\ref{M matrix}),
 then becomes
  \begin{align}
      M(t) &= 
         \mat{cc}{
		   M_{00} & M_{01} \\
		   M_{10} & M_{11} } =
     \mat{cc}{
      \Theta_0(0^+,t)+\Theta_0(1^-,t)  & \ 
      \Theta_\frac12(0^+,t)+\Theta_\frac12(1^-,t)  \\
      \Theta_0(\frac12^+,t)-\Theta_0(\frac12^-,t)  & \ 
      \Theta_\frac12(\frac12^+,t)-\Theta_\frac12(\frac12^-,t)  
		   }  \\
      & = \frac{1}{2} \ \mat{cc}
         { 
	   \frac{1}{1-t}  +\frac{1}{1-t}  & \ \
	   \frac{1}{1-t} + \frac{-1}{1+t}  \\[3mm]
	   \frac{1}{1-t}  -\frac{1}{1-t}  & \ \
	   \frac{1-2t}{1-t} - \frac{-1-2t}{1+t}  
	   }
      = \mat{cc}
         { \frac{1}{1-t} & \ \ \frac{-t}{1-t^2}  \\[3mm]
	   0   & \ \ \frac{1-t-t^2}{1-t^2} }
	   .
	 \end{align}
 The determinant is 
 $\ds \dMT(t)=\frac{1-t-t^2}{(1-t)^2 (1+t)}$ with one zero at 
 $t=2/(1+\sqrt{5})=1/\gamma$ 
 (implying $\htop=\log \gamma$)
 and no other zeroes in the unit disk.
 In the case of topological entropy it suffices in fact
 to consider the minor
 obtained from $M(t)$ by erasing the first line and the first column
 (see e.g.\ \cite[Appendix B]{RT}).

The corresponding Ruelle determinant $\dR(t)$ 
(or reciprocal zeta-function) is
calculated from the sequence \; $\# {\rm Fix} \whf^n$, $n\geq 1$,
the number of fixed points of the $n$-th iterate of the map 
(see e.g.\ \cite{BK}):
  \begin{align}
  \dR(t) 
     & = \exp \left(  - \sum_{n\geq 1} \frac{t^n}{n} 
      \# {\rm Fix} \whf^n \right) 
     = \exp \left(  - \sum_{n\geq 1} \frac{t^n}{n} 
       \tr T^n \right) \\
     &= \det \left( \bfone - t T \right) 
     = 1-t-t^2  ,
  \end{align}
where $T=\mat{cc}{1 & 1 \\ 1 & 0 }$ is the transition matrix
corresponding to 
$\whf \whI_1=\whI_1\cup\whI_2$ and
$\whf \whI_2=\whI_1$. The zeros
of this determinant are (see \cite{BK})
the peripheral eigenvalues of $L$ (counted with multiplicity).
     In accordance with our main theorem,
     $\dR(T)/\dMT(t) = (1-t)^2(1+t)$
     is analytic and without zeros in
     the unit disk.
We note that the precise form of $\dR(t)$ (but not the conclusion)
depends somewhat upon the choice of how to treat 
periodic orbits on the boundaries of the intervals.
In the above example, they may be omitted in the
sum, in which case one recovers
for $\dR(t)$
the (reciprocal of the) Lefschetz zeta function which is identical
to $\dMT(t)$ as  already shown in \cite{MT}.
\end{example}

Hans Henrik RUGH,      
University of Paris-Sud, 91405 Orsay Cedex, France.\\
e-mail: Hans-Henrik.Rugh@math.u-psud.fr

\end{document}